\documentclass[10pt]{article}
\input epsf

\begin{document}

\title{Hypercubes circumscribed in hyperspheres: a constant growth ratio for
volumes at large dimensions}
\small{\author{N. Shabot-Marcos$^a$ y A. Sandoval-Villalbazo$^a$  \\
$^a$ Departamento de F\'{\i}sica y Matem\'{a}ticas, Universidad Iberoamericana \\
Lomas de Santa Fe 01210 M\'{e}xico D.F., M\'{e}xico \\
E-Mail: alfredo.sandoval@uia.mx \\
E-Mail: nathanshabot@alumno.uia.mx }} \maketitle

\bigskip
\begin{abstract}
This note shows that an interesting property arises when
considering the relation between the hypersphere volumes at
dimensions $n+1$ and $n$, if the hyperspheres circumscribe unitary
hypercubes in $n+1$ and $n$ dimensions, respectively . In the
limit of large $n$, the ratio $\frac{V_{n+1}}{V_{n}}$ is constant
and equal to the irrational number $\sqrt{\frac{\pi e}{2}}$. Some
consequences of this result are explored.
\end{abstract}

It is well known that the  volume of a hypersphere  embedded in an
euclidean $n$-dimensional space decreases to zero for large $n$
\cite{uno} \cite{dos}. The value of $n$ in which the hypersphere
volume starts to decrease depends on the ratio $R$. On the other
hand, the corresponding volume of a unitary hypercube remains
constant in an euclidean $n$-dimensional space. We now consider
the problem of the circumscription of hypercubes in hyperspheres
in euclidean $n$-dimensional spaces.

The volume of a hypersphere in $n$ dimensions, for a given ratio
$R$, is given by:
\begin{equation}\label{one}
    V_{n}=\frac{2 \pi R^{n}}{n} \prod^{n-2}_{k=1}\frac{\sqrt{\pi}~\Gamma(
    \frac{k+1}{2})}{\Gamma(1+\frac{k}{2})},~n \geq3 .
\end{equation}
If the hypersphere circumscribes an unitary hypercube in the same
$n$ dimensional space, then, the hypersphere ratio is given by:
\begin{equation}\label{dos}
    R=\frac{\sqrt{n}}{2}
\end{equation}
So, the volumes $V_{n+1}$ and $V_{n}$ are related by the equation:
\begin{equation}\label{three}
    g_{n}=\frac{ V_{n+1}}{V_{n}}=\frac{n^{\frac{n}{2}-1}}{2}(n-1)^{\frac{3-n}{2}}
   \frac{\sqrt{\pi}~\Gamma(
    \frac{n-1}{2})}{\Gamma(\frac{n}{2})}
\end{equation}
Clearly, for large $n$:
\begin{equation}\label{four}
   \lim_{n\rightarrow \infty}g_{n}=\sqrt{\frac{\pi e}{2}}\simeq 2.0663656
\end{equation}
Fig. 1 shows, for continuous $n$, the corresponding plot of
Eq.(\ref{three})
\begin{figure}
\epsfxsize=3.4in \epsfysize=2.6in
\epsffile{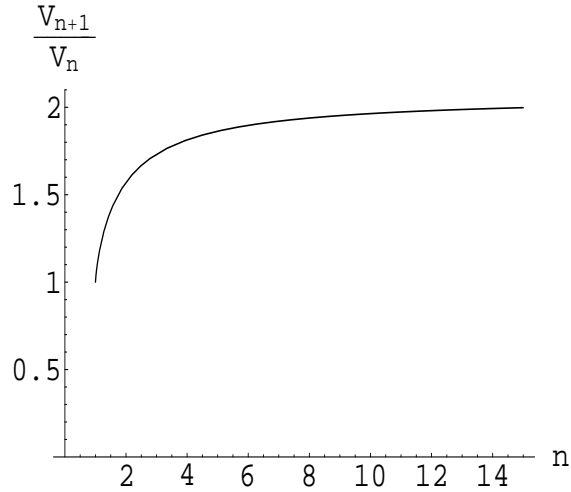}\vspace{0.5cm} \caption {The expression
$\frac{V_{n+1}}{V_{n}}$ for the case of hyperspheres of dimensions
$n+1$ and $n$ circumscribing unit hypercubes of the respective
dimensions.}
\end{figure}
\vspace{0.5cm}

To the authors' knowledge, this is the first time that the
irrational number $\sqrt{\frac{\pi e}{2}}$ appears in a
computation related to $n$-dimensional geometry. The result shows
that the circumscription requirement defines o family of
hyperspheres that possess constant volume growth as the number of
dimensions of its corresponding space tends to infinity. Further
properties of circumscription problems will be studied in future
work.

\end{document}